\documentclass[11pt]{article}
\usepackage{amsmath,amssymb,amsfonts}
\newtheorem{thm}{Theorem}[section]
\newtheorem{cor}[thm]{Corollary}
\newtheorem{lem}[thm]{Lemma}
\newtheorem{fct}[thm]{Fact}

\newtheorem{dfn}[thm]{Definition}
\newtheorem{pro}[thm]{Proposition}
\newtheorem{rem}[thm]{Remark}

=22cm\headheight=0cm\topskip=0cm\topmargin=0cm
\def\a{\bar {a}}\def\b{\bar {b}}
\def\x{\bar {x}}\def\y{\bar {y}}

\def\e{\mathbf {e}}

\def\proof{{\bf Proof\ }}

\begin{document}

\begin{center}
{\Large{\sc On some dynamical aspects of NIP theories 
}}

\vspace{15mm}

{\textbf{\bf Alireza Mofidi}} \vspace{2mm}

{\footnotesize Department of Mathematics and Computer Science,}
\\{\footnotesize Amirkabir University of Technology, Tehran, Iran}

and

{\footnotesize School of Mathematics, Institute for Research in Fundamental Sciences {\rm (IPM)}}
\\{\footnotesize{\it P.O. Box {\rm 19395-5746} Tehran, Iran.\\e-mail: mofidi@aut.ac.ir}}

\vspace{5mm}

\end{center}

\begin{abstract}
We study some dynamical aspects of the action of automorphisms in model theory
in particular in the presence of invariant measures.
We give some characterizations for NIP theories in terms of dynamics of automorphisms and invariant measures for example in terms of compact systems, entropy and measure algebras.
Moreover, we study the concept of symbolic representation for models.
Amongst the results, we give some characterizations for dividing lines and combinatorial configurations such as independence property, order property and strictly order property in terms of symbolic representations.

\end{abstract}

{\sc AMS subject classification 03C45, 37A05, 37B10, 37A35}

\bigskip

{\small {\sc Keywords}: NIP theories, measure preserving actions, dynamical systems, entropy.} 

\bigskip

\section{Introduction}

In his project for classification of first order theories, Shelah  
introduced several dividing lines in theories on the base of some combinatorial complexities. The machinery developed by him and others is called stability theory and nowadays its modern version, neo-stability, is used for analyzing classes beyond class of stable theories such as NIP theories. 
The class of NIP theories which are theories that no formula has independence property, is an important class containing stable theories and many studies in recent years were related to it.

On the other hand, many aspects of connection between definable groups (mostly in NIP context) with notions such as measure and forking are studied 
in for example \cite{HPP} and \cite{HP}. Also measures as technical tools were used in \cite{K1} for the first time in stability theory.
Let $G$ be a definable group acting on the topological space $S_G(M)$, space of types in $G$ with parameters from model $M$. This action gives a flow and some of its model theoretic and dynamical aspects are studied in \cite{N}. 
In that paper some definitions, notions and methodologies from topological dynamic 
are interpreted in the setting of model theory. 
Among them one can mention the notions of almost periodic types and Ellis semigroups.

In this paper we make some connections between some notions from dynamical systems in particular ergodic theory and symbolic dynamics to model theory. For example compact dynamical systems, entropy and symbolic representations are notions from ergodic theory and symbolic dynamics which are interacting with model theoretic notions in this paper.
We introduce some dynamical invariants and associate them to model theoretic objects (formulas, definable, models and theories, etc). Then we study 
certain model theoretic properties (for example stability theoretic properties) in terms of 
such associated notions.

One may equip a model or space of types with automorphism acting on them and also with invariant measure. Note that in the case of definable groups one may consider the action of the group on space of types. This structures forms some dynamical objects and one can study the model theoretic complexities such as stability configurations in terms of dynamical properties of such systems. 
We give several characterizations of NIP property in terms of dynamical and topological objects such as 
the measure algebras associated to models,  
abstract dynamical systems associated to the models when they are equipped with an automorphism and an invariant measure and also entropy. 
Among other invariants, we define and develop the concept of symbolic representation. In this case subflows of Bernoulli flow are associated structures.

Now we explain about the notion of symbolic representation.
It is a general philosophy in dynamical systems in particular symbolic dynamics that one associates some dynamical objects such as sub flows of 
Bernoulli flow $\mathcal{B}$ (namely $2^{\mathbb{Z}}$ 
equipped with the shift transformation) to general complicated dynamical systems. 
The goal is to study the properties of initial system via studying the properties of simpler associated ones.
This method was successfully used in several parts of dynamical system theory such as studying certain chaotic systems like horseshoe map by Smale.
Note that finding suitable symbolic representations of dynamical systems is a usual trends in dynamical systems.
In some part of this paper we more or less follow this general philosophy and associate some subflows of Bernoulli flows to model theoretic notions. 
As an example, for some instance $U$ of a formula $\phi(x,y)$, an automorphism $\sigma$ and $a \in M$ one may consider the subflows of $\mathcal{B}$ generated by sets of the form $\xi_{\sigma,U}(a):=\{n : \sigma^n(a) \in U\}$. We call sets of the form $\xi_{\sigma,U}(a)$ the symbolic image with respect to $\sigma$ and $U$.
Note that when one deals with the action of a single element of the acting group, 
then one can identify the orbit of a non periodic elements with $\mathbb{Z}$ and see the symbolic image as a subset of $2^{\mathbb{Z}}$.
Note that $2^{\mathbb{N}}$ is called the Cantor space. We use this name also for $2^{\mathbb{Z}}$.
We call functions such as $\xi_{\sigma,U}$ the symbolic representation functions. 
So symbolic representations associate sub flows of $\mathcal{B}$.
Studying NIP theories using the properties of such subflows is among interests. Also characterization of other combinatorial configurations such as order property strictly order property is also under consideration.

Organization of the paper is as follows.
In section \ref{sectionpreliminaries} we review some preliminaries.
In section \ref{Dynamicovermodels} we study NIP in terms of dynamic of action of automorphisms and invariant measures, in terms of entropy and also in terms of measure algebras.
We give some constructions for measures with certain properties and use them to give some characterizations for NIP properties in terms of compact dynamical systems and entropy. Also we study NIP theories from the point of view of measure algebras of models equipped with measures. In section \ref{sectionSymbolic view} we define the concept of symbolic representation for models and characterize NIP property, order property and strictly order property in terms of the notion of symbolic image. We also use these characterizations to obtain a theorem of Shelah which is connecting these notions together.

\section{Preliminaries and notations}\label{sectionpreliminaries}
In this section we review some definitions and facts which we need in the following sections.

\begin{dfn}\label{defnormalnumber}
\emph{Let $b$ be a natural number. A real number is called \textit{normal} in base $b$ if its sequence of digits in base $b$ is distributed uniformly where by distributed uniformly we mean that every single digits has the same density of appearance $1/b$, every tuple of digits has the same density of appearance $1/b^2$, etc. Also for the case $b=2$ and for every $p \in (0,1)$, we call a real number \textit{$p$-normal} (in base 2) if in its sequence of digits (in base 2), the density of appearance of any finite binary sequence $I$ of length $n$ is equal to $p^t(1-p)^{(n-t)}$ where $t$ is the number of appearance of $1$ in $I$.}
\end{dfn}

\begin{rem}\label{existenceofpnormalnumbers}
For every $p \in (0,1)$, the existence of $p$-normal numbers is proven in \cite{generalizedNormalNumber}.
\end{rem}

By a $\mathbb{Z}$-sequence we mean a sequence with indexes from $\mathbb{Z}$. In the following definition we associate a $\mathbb{Z}$-sequence to every non negative real number.
\begin{dfn}\label{sequencefromrealnmber}
\emph{For every non negative real number $u$ represented in the binary representation we define the \textit{$\mathbb{Z}$-sequence associated to $u$} to be a binary $\mathbb{Z}$-sequence $I_u$ with $I(i)=1$ if and only if $|i|$'th digit of $u$ is one for every $i \in \mathbb{Z}$. Note that one can see any binary $\mathbb{Z}$-sequence as a subset of $\mathbb{Z}$.}
\end{dfn}

\begin{dfn}
\emph{An \textit{abstract dynamical system} consists of a measure space $(X,\mathcal{A},\mu)$ equipped with a measure preserving map $\sigma: X \rightarrow X$. A \emph{topological dynamical system}, is a topological space, together with a continuous transformation, 
or more generally, a semigroup of continuous transformations of that space.} 
\end{dfn}

\begin{pro}\label{existenceinvmeasureoncmpttopspacehomo}
Let $X$ be a compact topological space and $\sigma:X \rightarrow X$ a homeomorphism. Then there exist a $\sigma$-invariant probability measure $\mu$ on Borel sigma-algebra of $X$.
\end{pro}

Now we review some notions from model theory in particular stability theory.
\begin{dfn}\label{defofOP} 
\emph{Let $\phi(x,y)$ be a formula. We say that $\phi$ has order property (OP) if there exists some model $M$ such that for every $n \in \mathbb{N}$, there exist $a_1,\ldots,a_n$ (with same arity of $x$) and $b_1,\ldots,b_n$ (with same arity of $y$) such that
$\phi(a_i,b_j) \Leftrightarrow i \leqslant j$.}
\end{dfn}

\begin{lem}\label{differentequivalencesofOP}
Let $T$ be a theory. Then the followings are equivalent.
\begin{enumerate}

\item{The formula $\phi(x,y)$ has OP.} 

\item{Let $I$ be an arbitrary linear order. Then there exists some model $M$ of $T$ and $\{a_i: i \in I\}$ and $\{b_i: i \in I\}$ sequences in $M$ such that 
$\phi(a_i,b_j) \Leftrightarrow i \leqslant j$.}\label{OPwitharbitraryOrder}

\item{Same as \ref{OPwitharbitraryOrder} with additional property that $\{a_i: i \in I\}$ is an indiscernible sequence.}
\end{enumerate}
\end{lem}

\begin{dfn}\label{defofIP}
\emph{Let $\phi(x,y)$ be a formula. We say that $\phi$ has independence property (IP) if there exists some model $M$ such that for every $n \in \mathbb{N}$, there exist $a_1,\ldots,a_n$ (with same arity of $x$)
such that witnesses IP of length $n$ i.e for every $J \subseteq \{1,\dots,n\}$, there exists some $b_J$ (with same arity of $y$) such that 
$\phi(a_i,b_J) \Leftrightarrow i \in J$.}
\end{dfn}

\begin{lem}\label{differentequivalencesofIP}
Let $T$ be a theory. Then the followings are equivalent.
\begin{enumerate}

\item{The formula $\phi(x,y)$ has IP. }

\item{Let $I$ be an arbitrary linear order. Then there exists some model $M$ of $T$ and sequences $\{a_i: i \in I\}$ in $M$ such that witnesses IP i.e
for every $J \subseteq I$, there exists some $b_J$ (with same arity as $y$) such that 
$\phi(a_i,b_J) \Leftrightarrow i \in J$.}\label{IPwitharbitrarysequence}

\item{Same as \ref{IPwitharbitrarysequence} with additional property that $\{a_i: i \in I\}$ is an indiscernible sequence.}\label{IPwitharbitraryindissequence}
\end{enumerate}
\end{lem}

\begin{dfn}\label{defofSOP} 
\emph{Let $\phi(x,y)$ be a formula. We say that $\phi$ has strict order property (SOP) if there exists some model $M$ and a sequence $(b_i)_{i < \omega}$ of 
$|y|$-tuples in $M$ such that $\phi(M,b_i) \subsetneqq \phi(M,b_{i+1})$ for every $i < \omega$. If $\phi(x,y)$ doesn't have SOP we call it NSOP. A theory is called a NSOP theory if every formula is NSOP.}
\end{dfn}

\begin{lem}\label{differentequivalencesofSOP}
\emph{
Let $T$ be a theory. Then the followings are equivalent.
\begin{enumerate}
\item{The formula $\phi(x,y)$ has SOP.}
\item{Let $I$ be an arbitrary linear order. Then there exists some model $M$ of $T$ and sequences $\{b_i: i \in I\}$ in $M$ such that 
$\phi(M,b_i) \subsetneqq \phi(M,b_{i+1})$ for every $i \in I$.}\label{SOPnestedsequence}
\item{Same as \ref{SOPnestedsequence} with additional condition that $\{b_i: i \in I\}$ is an indiscernible sequence.}
\end{enumerate}
}
\end{lem}

\noindent \textbf{Average types}
\bigskip

The notion of "average type" in NIP theories has many application. One can see 
\cite{Shelahclassification} for some of them.
We define the average type of a countable sequence.

\begin{dfn}\label{dfnavtypes}
\emph{Let $I =\{a_i : i < \omega \}$ be an indiscernible sequence and let $B$ be a subset of model. We define the average type of $I$ over parameter set $B$, $Av(I,B)$ as the set of all formula $\phi(x,b)$ with $b \in B$ such that $\{i : \neg \phi(a_i,b)\}$ is finite. 
} 
\end{dfn}
Note that under NIP assumption one has that $Av(I,A) \in S(A)$.

\bigskip

\noindent \textbf{Measures and NIP theories}

\bigskip

A (Keisler) measure on $M$ over parameter set $A$ is a finitely additive probability measure on  $Def_A(M)$, the definable sets of $M$ over parameter set $A$.

\begin{rem}\label{sigmaaddivitemeasureonmodel}
\emph{
The following is shown in \cite{K1}. Any Keisler measure on $M$ over parameter set $A$ can be uniquely extended to a countably additive measure on the sigma-algebra generated by the $Def(A)$ on $M$ (by Caratheodory’s theorem). 
Also a Keisler measure on $M$ over parameter set $A$ could be seen as a
regular Borel probability measure on the space of types $S(A)$.
}
\end{rem}

The following result of Keisler from \cite{K1} links Keisler measures and the notion NIP.

\begin{thm} \label{K-NIPKeislermeasures}
For a theory $T$ the following are equivalent:
\begin{enumerate}
\item{The theory $T$ has NIP.}\label{K-NIP}
\item{For every Keisler measure $\mu$ on $M$ and every formula $\phi(x,y)$, 
there is no infinite set of instances $\phi(x,a_i)$ such that the symmetric difference of the $\phi(x,a_i)$'s are bounded away from zero.}\label{K-bddaway0cbl}
\item{Same as \ref{K-bddaway0cbl} but with word 'finite' replaced by 'uncountable'.}\label{K-bddaway0finite}
\item{Every Keisler measure on $M$ has a countably generated measure algebra.}\label{K-cblgenmalg}
\end{enumerate}
\end{thm}

\section{Dynamical system on models and spaces of types}\label{Dynamicovermodels}
In this section we study some dynamical aspects of models and spaces of types when they are equipped with measures and automorphism. We first develop some tools for this.
Let $M$ be a model and $\sigma \in Aut(M)$. A subset $A \subseteq M$ is called $\sigma$-closed 
if $\sigma^n(A) \subset A$ for every $n \in \mathbb{Z}$. We let $\bar{A}=\bigcup_{i \in \mathbb{Z}}\sigma^i(A)$. Obviously $\bar{A}$ is closed under action of $A$ and if the set $A$  is small ($|A| < |M|$) then $\bar{A}$ is small too.

\begin{rem}
\emph{
Let $M$ be a model, $\sigma \in Aut(M)$ and $A \subseteq M$ be $\sigma$-closed.
Then $\sigma$ is measurable with respect to the Borel sigma-algebra generated by $Def(A)$. 
}
\end{rem}

Let $M$ be a model, $\sigma \in Aut(M)$ and $A \subseteq M$ be $\sigma$-closed. 
Then $\sigma$ can be seen an homeomorphism acting on space of types $S(A)$
where $S(A)$ is a compact Hausdorff space with logic topology $\tau$.
Hence $(S(A),\tau,\sigma)$ is a topological dynamical system. Moreover for every $\sigma$ invariant measure $\mu$ one can see $\mu$ as a Borel probability measure on $S(M)$ (by Remark \ref{sigmaaddivitemeasureonmodel}) and so $(S(M),\mu,\sigma)$ is a abstract dynamical system.
Note that the existence of $\sigma$-invariant measures on $S(A)$ is guaranteed 
by Proposition \ref{existenceinvmeasureoncmpttopspacehomo}.

Now we talk about a method for constructing measures with specific properties.
We will use it later in the following sections.
Let $\mathcal{F}$ be an  
ultrafilter extending Frechet filter on $\mathcal{P}(\mathbb{N})$. We use the notation $\lim^{\mathcal{F}}$ for the notion of limit of sequences with respect to the ultrafilter $\mathcal{F}$. More precisely for every real sequence $(a_i)$, $\lim^{\mathcal{F}}(a_i)$ is the unique element $c$ such that for every $\epsilon >0$ we have that $\{i: |c-a_i| <\epsilon\} \in \mathcal{F}$. Note that since $\mathcal{F}$ is extending Frechet filter, for any convergent sequence $J=(a_i)$ in 
$\mathbb{R}$, $\lim(J)$ and $\lim^\mathcal{F}(J)$ will be the same.

\begin{dfn}
\emph{For every $I \subseteq \mathbb{Z}$ define $dns(I):=\lim_{n \rightarrow \infty} \frac{I \cap [-n,n]}{n}$ provided that limit exists.
Let $\mathcal{F}$ be an  
ultrafilter extending Frechet filter on $\mathcal{P}(\mathbb{N})$.
Define $dns_{\mathcal{F}}(I):=\lim_{n \rightarrow \infty}^\mathcal{F} \frac{I \cap [-n,n]}{n}$
where by $\lim_{n \rightarrow \infty}^\mathcal{F}$ we mean the limit with respect to the ultrafilter $\mathcal{F}$.}
\end{dfn}

\begin{dfn}\label{defdnsWI}
\emph{Let $I$ be a $\mathbb{Z}$-sequence and $W$ be a finite sequence. We define 
$\langle W,I \rangle:=\{i:  W(j) =I(i+j) \ for \ all \  j=1,\ldots, |W| \}$.
So $dns_{\mathcal{F}}(\langle W,I \rangle)$ roughly speaking is the density of appearance of $W$ in $I$.}
\end{dfn}

\begin{rem}\label{capI-j=<W,I>}
Let $I$ and $W$ be as Definition \ref{defdnsWI}. Then we have that
$$\bigcap_{1 \leqslant j \leqslant |W|} (I-j)^{W(i)}=\langle W,I \rangle.$$
\end{rem}

\begin{dfn}\label{defsymbimagsingleelement}
\emph{Let $X$ be a set, $B \subseteq X$, $a \in X$ and $\sigma:X \rightarrow X$ be a map.
We define $\xi_{\sigma,B}(a):=\{i \in \mathbb{Z}: \sigma^i(a) \in B\}$.}
\end{dfn}

\begin{dfn}\label{dfnlimitfrequencymeasures}
\emph{Let $\mathcal{F}$ be an ultrafilter on $P(\mathbb{N})$,
$\a=(a_i: i \in \mathbb{Z})$ be a sequence in $X$ and $\textbf{B}$ be a Boolean algebra of subsets of $X$. 
We define a function $\mu^{\a}_\mathcal{F}$ on $\textbf{B}$ in the following way.
For each $B \in \textbf{B}$ let 
$\mu^{\a}_\mathcal{F}(B)=dns_{\mathcal{F}}(\{i:a_i \in B\})$.
In particular $(a_i: i \in \mathbb{Z})$ can be the orbit of some $a \in X$ under iteration of some map $\sigma$ from $X$ to itself. So
for each $B \in \textbf{B}$ we have 
$\mu^{\a}_\mathcal{F}(B)=dns_{\mathcal{F}}(\xi_{\sigma,B}(a))$.
In this case we use the notation $\mu^{\sigma,a}_\mathcal{F}$ instead of $\mu^{\a}_\mathcal{F}$.
We call $\mu^{\a}_\mathcal{F}$ and $\mu^{\sigma,a}_\mathcal{F}$ the \textit{limit frequency measures} on $\textbf{B}$ with respect to the sequence $\a$.
}
\end{dfn}

Note that $\mu^{\sigma,a}_\mathcal{F}$ defined above is a finitely additive measure on $\textbf{B}$ and is $\sigma$-invariant.

\begin{rem}
\emph{
Let $M$ be a saturated enough model. Then the average types (defined in \ref{dfnavtypes}) are special cases of limit frequency measures.
}
\end{rem}
\proof
Let $I=(c_i)_{i \in \mathbb{Z}}$ be an indiscernible sequence over $A \subseteq M$.
By saturation, there exists some automorphism $\sigma \in Aut(M/A)$ and $a \in M$ such that 
$I$ is some part of $\sigma$-orbit of $a$. Let $\mathcal{F}$ be an extension of Frechet filter to an ultrafilter.
Now consider $\mu_{\mathcal{F}}^{\sigma,a}$ on $\textbf{B}=Def(M)$. Then the set of those definables $D$ with $\mu_{\mathcal{F}}^{\sigma,a}(D)=1$ are those containing co-finite part of 
$I$ and therefore belong to $Av(I,A)$.
$\ \ \square$

\subsection{NIP theories}
In this subsection we study effects of the property NIP on topological space and 
dynamical systems which we associate to models.

\subsubsection{NIP and compact abstract dynamical systems}\label{NIPcptADS}
In this section we will characterize NIP theories in terms of compact (almost periodic) dynamical systems in ergodic theory.

\begin{dfn} \label{defCptADS}
\emph{We say that the abstract dynamical system $(X,\mathcal{A},\mu,\sigma)$ is a compact system (almost periodic system) if for every $\epsilon >0$ and $A \in \mathcal{A}$, there exists some $n \in \mathbb{N}$ such that $\mu(\sigma^n(A) \triangle A) <\epsilon$.}
\end{dfn}
Note that compact or almost periodic dynamical systems are closest ones to periodic systems in the hierarchy of dynamical systems.

\begin{dfn}
\emph{A subset $I \subseteq \mathbb{Z}$ 
 is called $\epsilon$-wide if for every $n \in \mathbb{Z}$ 
we have that $$dns(I)-dns(I \cap I_n)\geqslant \epsilon$$ 
where $I_n$ is obtained from $I$ by $n$ shifts.}
\end{dfn}

\begin{rem}\label{ubsetN}
\emph{
Let $u$ be a normal real number (in base 2) and $I_u$ be 
the $\mathbb{Z}$-sequence associated to $u$
(defined in Definition \ref{sequencefromrealnmber}).
Then $I_u$ is $\epsilon$-wide for $\epsilon=\frac{1}{4}$.
}
\end{rem}

\begin{dfn}
\emph{Let $\phi(x,y)$ be a formula.
\begin{enumerate}
\item{$\phi(x,y)$ is called strong wide in a model $M$ of a theory $T$ if there exist $a \in M$, $c \in M$, $\sigma \in Aut(M)$ and $\epsilon>0$ such that $\xi_{\sigma,U}(a)$ 
is a $\epsilon$-wide subset of $\mathbb{Z}$ where $U=\phi(x,c)$.
A theory $T$ is called strong wide if there exists some formula $\phi(x,y)$ such that is strong wide in some model of $T$.
}
\item{$\phi(x,y)$ is called weak wide in a model $M$ of a theory $T$ if there exist two sequences $(a_n)_{n \in \mathbb{Z}} \in M$, $(c_n)_{n \in \mathbb{Z}} \in M^t$ and some $\epsilon>0$ such that letting $U_n=\{a_i: \phi(a_i,c_n)\}$ for each $n \in \mathbb{Z}$, we have that
$U_n=\{a_{i+n}:a_i \in U_0 \}$ for each $n \in \mathbb{Z}$ and $\{i: a_i \in U_0\}$ is a $\epsilon$-wide subset of $\mathbb{Z}$. A theory $T$ is called weak wide if there exists some formula $\phi(x,y)$ such that is weak wide in some model of $T$.
}
\end{enumerate}
}
\end{dfn}

The following statement gives a characterization of the NIP in terms of compact dynamical systems.
In the following theorem we consider our measures to be countably additive.
Note that by Remark \ref{sigmaaddivitemeasureonmodel} a measure on $Def(M)$ can be extended to a countably additive Borel measure on the sigma-algebra generated by $Def(M)$.
\begin{thm}\label{NIPvsCPTADS}
Let $M$ be a saturated model of theory $T$. 
Then the followings are equivalent:
\begin{enumerate}

\item{$T$ has IP.}\label{TisNIP}

\item{There is some nontrivial automorphism $\sigma$ of $M$, some $A \subseteq M$ and some $\sigma$-invariant measure $\mu$ such that $(M, \mathcal{A}, \mu,\sigma)$ is not a compact dynamical system where $\mathcal{A}$ is the sigma-algebra generated by $Def(A)$. 
}
\label{everymeasureCPTDS}

\item{Same as \ref{everymeasureCPTDS} except that the measure $\mu$ is a global measure.}
\label{everyglobalmeasureCPTDS}

\item{ $T$ is weak wide.}
\label{notweakwide}

\item{ $T$ is strong wide.}
\label{notstrongwide}

\end{enumerate}

\end{thm}
\proof
$(\ref{TisNIP} \Rightarrow \ref{notstrongwide} )$
Assume 
that $T$ has IP. Then by a theorem of Shelah (see for example Theorem 12.18 of \cite{Poizatmodeltheory}),  
there exists some formula $\phi(x,y)$ which witnesses IP with $|x|=1$. Assume that $|y|=t$.
By part \ref{IPwitharbitraryindissequence} of the Lemma \ref{differentequivalencesofIP} there exists an indiscernible sequence $\{a_i: i \in \mathbb{Z}\}$ (where $a_i$'s are of arrity 1) and a sequence $\{a_i: i \in \mathbb{Z}\}$ witnessing IP for $\phi$. Since $M$ is saturated and $\{a_i: i \in \mathbb{Z}\}$ are indiscernible, there exists some $a \in M$ and some $\sigma \in Aut(M)$ such that $\sigma^i(a)=a_i$ for every $i \in \mathbb{Z}$.
Now using Remark \ref{ubsetN}, we find some arbitrary $\epsilon$-wide subset of $\mathbb{Z}$, 
say $I $. Since $\phi$ has IP property and $M$ is saturated, there exists some $c \in M^t$ such that $\phi(a_i,c) \Leftrightarrow i \in I$. So $\xi_{\sigma,U}(a)=I$ where $U=\phi(M,c)$.
Hence $M$ is strong wide. Therefore $T$ is strong wide. 

\vspace{1.5mm}

$(\ref{notstrongwide} \Rightarrow \ref{everyglobalmeasureCPTDS})$
Assume 
that $M$ is a model of $T$ which is strong wide witnesses by a formula $\phi(x,y)$, $a \in M$, $c \in M^t$, $\epsilon>0$ and automorphism $\sigma$.
So $\xi_{\sigma,U}(a)$
is a $\epsilon$-wide subset of $\mathbb{Z}$ where $U=\phi(M,c)$. 
Let $\mathcal{F}$ be an 
ultrafilter extending Frechet filter over $\mathbb{N}$ and let 
$\mu^{\sigma,a}_{\mathcal{F}}$
be the corresponding limit frequency measure defined in the Definition \ref{dfnlimitfrequencymeasures} on $\mathcal{A}$, the sigma-algebra generated by $Def(M)$.
We use notation $\mu$ instead of $\mu^{\sigma,a}_{\mathcal{F}}$.
Clearly $\mu$ is $\sigma$-invariant.
Now we show that $(M,\mathcal{A},\mu,\sigma)$ is not a compact dynamical system.
For every $n \in \mathbb{Z}$, let $U_n=\sigma^n(U)$, $I=\xi_{\sigma,U}(a)$ and $J=\xi_{\sigma,U\cap U_n}(a)$. 
We have $\mu(U)=\mu(U_n)$ and $J=I \cap I_n$.
So we have that 
$$\mu (U \triangle U_n)=\mu(U \setminus U_n)+\mu(U_n\setminus U)=2\mu(U)-2\mu(U \cap U_n)=2dns_{\mathcal{F}}(I)-2dns_{\mathcal{F}}(J) \geqslant 2\epsilon.$$
Therefore $(M,\mathcal{A},\mu,\sigma)$ is not a compact dynamical system.

\vspace{1.5mm}

$(\ref{everyglobalmeasureCPTDS} \Rightarrow \ref{everymeasureCPTDS})$ and 
$(\ref{notstrongwide} \Rightarrow \ref{notweakwide})$
are Obvious.

\vspace{1.5mm}

$(\ref{everymeasureCPTDS} \Rightarrow \ref{TisNIP})$
We prove by contradiction. Assume that $T$ is NIP. Let $B=\phi(M,b)$ be an arbitrary instance of a formula and let $\epsilon >0$ be arbitrary. Note that for every $n$ we have
$\sigma^n(B)=\phi(M,\sigma^n(b))$. By using Theorem \ref{K-NIPKeislermeasures} we find $N_{\epsilon}$ such that
$\forall m > n > N_{\epsilon}, \ \mu(\sigma^n(B) \triangle \sigma^m(B))< \epsilon$.
Since $\mu$ is $\sigma$-invariant, we have 
$\mu(B \triangle \sigma^{m-n}(B))=\mu(\sigma^n(B) \triangle \sigma^m(B))< \epsilon$.
Since the measures of definable sets determines measures of all measurable's, then $\sigma$ is almost periodic and makes the system a compact dynamical system which is a contradiction. 

One can use  
the Proposition $\ref{equivalentofNIPintermMeasureAlgebra}$
to give another proof by contradiction for this part. Assume that $\phi$ has NIP. Hence for every measure, the corresponding measure algebra $\mathcal{M}$ is precompact. Moreover $\sigma$ induces an isometry on $\mathcal{M}$. We denote the element of $\mathcal{M}$ corresponding to $\sigma^{i}(B)$ by $b_i$. Also denote the metric on $\mathcal{M}$ obtained from $\mu$ by $d_{\mu}$. Since $\mathcal{M}$ is precompact the sequence $b_i$ has a Cauchy subsequence. Thus the are $m,n$ such that $d_{\mu}(b_n,b_m)<\epsilon$. Hence 
$$\mu(B \triangle \sigma^{m-n}(B))=d_{\mu}(b_0,b_{m-n})=d_{\mu}(b_n,b_m) < \epsilon.$$
This shows that system is a compact dynamical system which is a contradiction.

\vspace{1.5mm}

$(\ref{notweakwide} \Rightarrow \ref{TisNIP})$
Assume 
that $M$ is a model of $T$ which is weak wide witnesses by a formula $\phi(x,y)$, $\epsilon>0$ and sequences $a_n \in M$, $c_n \in M^t$.
Recall that $U_n=\{a_i: \phi(a_i,c_n)\}$ for every $n$. Also $U_0$ is $\epsilon$-wide and $U_n=\{a_{i+n}:a_i \in U_0 \}$ for each $n \in \mathbb{Z}$.
Let $I_n=\{i: a_i \in U_n\}$. Clearly $I_n$ can be obtained from $I_0$ by $n$ shifts.
Also let $\a$ denotes the sequence $(a_i)_{\mathbb{Z}}$ and let $\mathcal{F}$ be a 
ultrafilter extending the Frechet filter over $\mathbb{N}$ and let 
$\mu:=\mu^{\a}_\mathcal{F}$ be limit frequency measure defined in the Definition \ref{dfnlimitfrequencymeasures} with respect to the sequence $\bar{a}$. 
Let $D_n:=\phi(M,c_n)$ for every $n$. One can see that for every $m,n$ we have that
$\mu(D_n)=\mu(D_m)$. So we have that
$$\mu(D_n \bigtriangleup D_m)=\mu(D_n \setminus D_m)+\mu(D_m \setminus D_n)=
2\mu(D_n)-2\mu(D_n \cap D_m)$$
$$=2dns_{\mathcal{F}}(I_n)-2dns_{\mathcal{F}}(I_n \cap I_m) \geqslant 2\epsilon.$$
Now using the theorem \ref{K-NIPKeislermeasures}, $\phi$ has IP.
Note that we have used the facts that $dns_{\mathcal{F}}(I_n)=dns_{\mathcal{F}}(I_m)$ and $dns_{\mathcal{F}}(I_n \cap I_m)=dns_{\mathcal{F}}(I_0 \cap I_{m-n})$.
$\ \ \square$

\begin{example}
Let $M$ be a model of the theory of random graphs. By quantifier elimination every definable in $M$ could be seen as disjunctions of formulas of the form
$$D=R(x,a_1) \wedge \ldots \wedge R(x,a_n) \wedge \neg R(x,b_1) \ldots \wedge \neg R(x,b_m)$$
for elements $a_i$'s and $b_j$'s in $M$. 
We denote by $Def(M)$ the set of all definable's in $M$.
A Bernoulli measure with parameter $p \in (0,1)$ on $Def(M)$ is defined with the following valuation for every $D$.
$$\mu(D)=p^n(1-p)^m.$$
One sees that $\mu$ is an $Aut(M)$-invariant measure on $M$.
Also it is not hard to see that for any $\sigma \in Aut(M)$, $(M,Def(M),\mu,\sigma)$ forms a non-compact dynamical system
and so theorem \ref{NIPvsCPTADS} verifies existence of IP in the theory of random graphs.
Note that in here we use the notion of dynamical system on a finitely additive measure in a similar way of countably additive ones.

\end{example}

\subsubsection{NIP and entropy}\label{NIPandentropy}
In this part we give a characterization of NIP in terms of measure theoretic entropy of automorphism of models. Roughly speaking we show that NIP is equivalent to entropy zero on every automorphisms. We first briefly review some notions of measure theoretic entropy. For a more extensive description one can see for example \cite{EntropyinDS}.

Let $(M,\mathcal{A},\mu)$ be a measure space and $\sigma$ be a measurable map on $M$ such that $\mu$ is $\sigma$-invariant. For every measurable partition $P=\{P_1,\ldots,P_n\}$ we define the measure theoretic entropy of $P$ with respect to $\sigma$ and $\mu$ with $h(\sigma,P)=\lim_{n \rightarrow \infty} \frac{1}{n} H(\bigvee_{i=0}^{n-1}\sigma^{-i}(P))$
where $H(\bigvee_{i=0}^{n-1}\sigma^{-i}(P))=\sum p_i \log p_i$ where $p_i$'s are the measures of atoms of the partition $\bigvee_{i=0}^{n-1}\sigma^{-i}(P)$.
Now entropy of $\sigma$ is defined as $h(\sigma)=\sup_{P} h(\sigma,P)$.

\vspace{1.5mm}

We will need the following lemmas later.
\begin{lem}\label{boundonHforpartition}
Let $P$ be a partition with $k$ atoms. Then we have that $H(P) \leqslant \log k$.
\end{lem}

\begin{lem}\label{ent0impliesent0ongenerated}
Assume that $X$ is a set (possibly equipped with additional structures), $\sigma$ is an automorphism of $X$ and $\textbf{A} \subseteq \mathcal{P}(X)$. Also let $\mathcal{A}$ be the sigma algebra generated by $\textbf{A}$ and let $\mu$ be a $\sigma$-invariant measure on $\mathcal{A}$. Assume that for every finite partition $P$ of elements of $\textbf{A}$
we have that $h(\sigma,P)=0$. Then we have that $h(\sigma)=0$.
\end{lem}
\proof
Let $\Omega$ be the set of all countable ordinals. 
We have that $\mathcal{A}=\bigcup_{\alpha \in \Omega} \mathcal{E}_{\alpha}$ where for limit ordinals $\alpha$ we define $\mathcal{E}_{\alpha}:=\bigcup_{\beta < \Omega} \mathcal{E}_{\beta}$ and for successor ordinal $\alpha=\beta+1$ we define 
$\mathcal{E}_{\alpha}$ to be the collection of sets of the form of countable unions of elements of $\mathcal{E}_{\beta}$ or complements of such. 
By a transfinite induction on $\Omega$ and also a using the continuity of entropy function on space of partitions (see for example fact 1.7.9 of \cite{EntropyinDS}) the proof will be obtained.
$\square$

\bigskip

By using the Remark \ref{sigmaaddivitemeasureonmodel}, in the following theorem by a measure we mean the unique countably additive extension to $\langle Def(M) \rangle$ (the sigma-algebra generated by $Def(M)$) of the measure on $Def(M)$ .

\begin{thm}\label{NIPiffazeforaut}
Let $M$ be a saturated model of the theory $T$. Then $T$ is NIP if and only if for every $\sigma \in Aut(M)$ and every $\sigma$-invariant measure $\mu$ on $\langle Def(M) \rangle$, we have that 
$h_{\mu}(\sigma)=0$.
\end{thm}
\proof
$(\Rightarrow)$
By using \ref{ent0impliesent0ongenerated} it is enough to show that every partition of definable sets have zero entropy.
We first show that NIP implies a.z.e (i.e entropy of every $\sigma \in Aut(M)$ w.r.t every $\sigma$-invariant measure is zero).
Let $P=\{\phi_1(x,a_1),\ldots,\phi_r(x,a_r)\}$ be a definable partition of $M$. We show that $h(\sigma,P)=0$.
We have that 
$$h(\sigma,P)=\lim_{n \rightarrow \infty} \frac{1}{n} H(\bigvee_{i=0}^{n-1}\sigma^{-i}(P)).$$
But $\bigvee_{i=0}^{n-1}\sigma^{-i}(P)$ is exactly same as $S_{\Delta}(A)$ where 
$\Delta=\{\phi_1,\ldots,\phi_r\}$ and
$$A=\{a_1,\ldots,a_r,\sigma^{-1}(a_1), \ldots, \sigma^{-1}(a_r), \ldots, \sigma^{-n}(a_1), \ldots, \sigma^{-n}(a_r)\}.$$
Using a facts about NIP theories 
there exists some constant real $t$ such that 
$S_{\Delta}(A)\leqslant |A|^t$.
So using Lemma \ref{boundonHforpartition} for every $n$ we have that 
$$H(\bigvee_{i=0}^{n-1}\sigma^{-i}(P))\leqslant log|A|^t=t.log (r(n+1)).$$ 
Hence 
$$h(\sigma,P)\leqslant \lim_{n \rightarrow \infty} \frac{1}{n}.t.log (r(n+1))=0$$
and the proof is complete.
\bigskip

$(\Leftarrow)$
We prove by contradiction that entropy zero for all automorphisms with respect to invariant measures implies NIP. So assume that $T$ has IP.
Similar to the proof of part $(\ref{TisNIP} \Rightarrow \ref{notstrongwide})$ of the Theorem \ref{NIPvsCPTADS},
we can find a formula $\phi(x,y)$ ($|x|=1$ and $|y|=t$) witnessing IP in $M$
by a sequence $\mathcal{O}(a)=\{a_n\}$ where $a_n=\sigma^n(a)$ for each $n \in \mathbb{Z}$ where $\sigma$ is an automorphism of $M$.
Also we let $I \subseteq \mathcal{O}(a)$ be the $\frac{1}{4}$-wide subset of $\mathbb{Z}$ obtained from a normal number by the method explained in Remark \ref{ubsetN} and 
$c \in M^t$ be such that $\phi(a_n,c) \Leftrightarrow n \in I$.
Using the same method of proof of part $(\ref{notstrongwide} \Rightarrow \ref{everyglobalmeasureCPTDS})$ 
of the same theorem, we construct a limit frequency measure $\mu$ obtained from $\mathcal{O}(a), \sigma$ and $I$. 
We claim that $(M,\mu,\sigma)$ has positive entropy. For that we prove that the measurable partition  
$P=\{P_0,P_1\}$ has positive entropy, where $P_1=\phi(M,c)$ and $P_0=P_1^c=\neg \phi(M,c)$. Note that 
$\mu(P_0)=\mu(P_1)=dns(I)=\frac{1}{2}$.
For each $n$, the partition 
$\bigvee_{i=0}^{n-1}\sigma^{-i}(P)$
consists of sets of the form $\bigcap_{i=0}^{n-1} \sigma^{-i}(P_{v_i})$ for $v=(v_1,\ldots,v_{n-1}) \in \{0,1\}^{n-1}$.
But since $I$ was corresponding to a normal number, for every $v=(v_1,\ldots,v_{n-1}) \in \{0,1\}^{n-1}$ we have that 
$$\mu (\bigcap_{i=0}^{n-1} \sigma^{-i}(P_{v_i}))=dns(\bigcap_{i=0}^{n-1} (I^{v_i}-i))=\frac{1}{2^n}$$ 
where $I^1=I$ and $I^0$ is obtained by replacing $0$ and $1$ in $I$ and also by $I^{v_i}-i$ we mean $i$ times shift of $I^{v_i}$. So we have that 
$$H(\bigvee_{i=0}^{n-1}\sigma^{-i}(P))=-2^n.\frac{1}{2^n}.log(\frac{1}{2^n}).$$
Therefore $h(\sigma,P)=1$ which is a contradiction.
Now the proof is complete.
$\ \ \ \square$

\vspace{2mm}

Now we want to give an analogue of Theorem
\ref{NIPiffazeforaut} for definable groups in which the action of automorphisms are replaces by action of group.

\begin{rem}\label{inNIPdefinableGroupsareAZEwithgrouptranslation}
\emph{
Let $G$ be a definable 
group in a model $M$ of a NIP theory. Also let $\mu$ be an invariant measure on $G$. Then the entropy of action of every element of $G$ 
(by translation) with respect to $\mu$ is zero.
}
\end{rem}
\proof
Let $\tau$ be action of $g \in G$ on $G$ by translation (i.e $\tau(U)=g.U$ for every definable $U$). 
Let $\mathcal{U}=\{\phi^G_1(x,a_1),\ldots,\phi^G_r(x,a_r)\}$ be a definable partition of $G$ where $a_i \in M$.
Now rewrite the proof of \ref{NIPiffazeforaut}(i) just by replacing $\sigma$ with $\tau$.
$\ \ \square$

\vspace{2mm}

Now we look at to an example of a basic and fundamental object in dynamical systems from a model theoretic point of view.

\begin{example}\label{S1example}
The group $\mathbb{T}=(S^1,+)$ is a definable group in O-minimal structure $(\mathbb{R},+,.,0,1)$ which it's theory is NIP.
the action of every element of $\mathbb{T}$ on $\mathbb{T}$ by group operation, gives us a rotation on $\mathbb{T}$. Dynamic of homeomorphism of
$\mathbb{T}$ are studied extensively and completely characterized by Poincare. It is also known in 
that entropy of rotations on $\mathbb{T}$ is zero.
\end{example}

\subsubsection{NIP and Measure algebra}
We state the following fact and use it for describing some topological properties of associated measure algebras to models of NIP theories. We will need the following lemma later.
\begin{lem}\label{everyseqhasCauchysubseq}
Let $(X,d)$ be a metric space with the property that for each sequence $\{a_i\}$ of elements of $X$, the set of values $\{d(a_i,a_j):i \neq j\}$ is not bounded away from zero. Then every sequence in $X$ has a Cauchy subsequence.
\end{lem}

\begin{dfn}
\emph{
Let $\Lambda=(X,\mathcal{F},\mu)$ and
$\Gamma=(Y,\mathcal{G},\nu)$  be two probability spaces. 
\begin{enumerate}
\item{By an \textit{isomorphism} between $\Lambda$ and $\Gamma$ we mean an invertible map 
$f : X \rightarrow Y$  such that $f$  and $f^{-1}$ are measurable and measure preserving. In that case,, we call $\Lambda$ and $\Gamma$ \textit{isomorphic}.
Also we call $\Lambda$ and $\Gamma$ \textit{isomorphic mod 0} (or \textit{point isomorphic}) 
 if there exist null sets $X_1 \subset X$ and $Y_1 \subset Y$  such that the probability spaces $X \setminus X_1$ and $Y \setminus Y_1$ are isomorphic. 
}
\item{The probability space $(X,\mathcal{F},\mu)$ is called \textit{standard} (or \textit{Lebesgue spaces}) if it is isomorphic mod 0 either to an interval with Lebesgue measure, or a finite or countable set of atoms, or a combination (disjoint union) of both.
}
\item{By a \textit{set isomorphism} between $\Lambda$ and $\Gamma$ we mean a measure preserving isomorphism between their measure algebras. If such map exists we call $\Lambda$ and $\Gamma$ set isomorphic. 
}
\end{enumerate}
}
\end{dfn}
The theory of standard probability spaces (or Lebesgue spaces) was initiated by von Neumann and developed by Rokhlin. 
We denote the measure algebra of equivalence classes of Lebesgue measurable subsets of the interval $[0,1]$ with $\Omega$.
The following important statement (which is usually called "the isomorphism theorem") classifies the measure algebras in terms of set isomorphism.

\begin{fct}\label{isomorphismthmsmeasurespacesalgebrasHalmosVonNeumanthm1} 
A probability space is set isomorphic to the unit interval if and only if it is separable and non atomic. (see Theorem $1$ in \cite{HalmosVonNeuman})
\end{fct}

The following statement is a characterization of NIP in terms of measure algebras corresponding to models. We show the connection between NIP theories and standard probability spaces.

\begin{pro}\label{equivalentofNIPintermMeasureAlgebra}
The followings are equivalent.
\begin{enumerate}
\item{The theory $T$ has NIP.}\label{equivalentofNIPintermMeasureAlgebraThasNIP}
\item{For every model $M$ and every measure $\mu$ on $Def(M)$, the measure algebra corresponding to $(M,Def(M),\mu)$ is a precompact topological space.}\label{equivalentofNIPintermMeasureAlgebraNIPiffMALGisprecpt}
\item{For every model $M$, every $A \subset M$ and every non-atomic measure $\mu$ on $\mathcal{A}$, the sigma-algebra generated by $Def(A)$,  the measure space $(M, \mathcal{A}, \mu)$ is set isomorphic to some Lebesgue space (equivalently set isomorphic to the unit interval).}\label{equivalentofNIPintermMeasureAlgebraNIPsetisounitint}
\end{enumerate}
\end{pro}
\proof
$(\ref{equivalentofNIPintermMeasureAlgebraThasNIP} \Leftrightarrow \ref{equivalentofNIPintermMeasureAlgebraNIPiffMALGisprecpt})$
By using Theorem \ref{K-NIPKeislermeasures} and Lemma \ref{everyseqhasCauchysubseq}. 

\vspace{1.5mm}

$(\ref{equivalentofNIPintermMeasureAlgebraThasNIP} \Rightarrow \ref{equivalentofNIPintermMeasureAlgebraNIPsetisounitint})$
Assume that $T$ has NIP and $\mu$ is a measure on $Def(A)$ on model $M$. By part \ref{K-cblgenmalg} of the Theorem \ref{K-NIPKeislermeasures}, the measure algebra corresponding to $\mu$ is countably generated and therefore is separable. Now by using 
of the Fact \ref{isomorphismthmsmeasurespacesalgebrasHalmosVonNeumanthm1},
it is set isomorphic to a Lebesgue space.

\vspace{1.5mm}

$(\ref{equivalentofNIPintermMeasureAlgebraNIPsetisounitint} \Rightarrow \ref{equivalentofNIPintermMeasureAlgebraThasNIP})$
Let $\mu$ be an arbitrary measure on $Def(A)$. If $\mu$ is non-atomic then by assumption 
$(M, \mathcal{A}, \mu)$ is set isomorphic to some (any) Lebesgue space. So its measure algebra $\mathcal{M}$ is isomorphic to measure algebra of a Lebesgue space, say unit interval.
Therefore $\mathcal{M}$ has a countable generator. 
If $\mu$ has nontrivial atomic and non-atomic parts, 
one can ignore the measure of the atomic part, then normalize the non-atomic part and use the above argument to show that measure algebra is countably generated. Completely atomic measures also obviously have countably generated measure algebras. Now by using part $\ref{K-cblgenmalg}$ of Theorem $\ref{K-NIPKeislermeasures}$ the theory $T$ has NIP. 
$\ \square$

\section{A symbolic dynamical view to NIP theories}\label{sectionSymbolic view}
In this section we give characterizations for some stability theoretic classes such as stable, NIP and NSOP 
in terms of symbolic representation of single automorphisms. We also obtain a theorem of Shelah using these characterizations.
Symbolic dynamic has several application in analyzing dynamical systems in ergodic theory, topological dynamic and even algebraic dynamics.
As an example, in algebraic dynamical systems one usually deals with some $f$ from an algebraic variety $f$ to itself. In this situation one works with the orbit of the points, namely $\{f^n(a), n \in \mathbb{N}\}$ for every point $a$ of the variety. Also an important notion is the the notion of return sets of a point $a$ corresponding a set $Y$, namely 
$\{n \in \mathbb{Z}, f^n(a) \in Y\}$.
Using the Definition \ref{defsymbimagsingleelement}, one sees that this notion is exactly same as $\xi_{f,Y}(a)$.

Now we review some notions from dynamical systems.
Let $\mathcal{B}=(2^{\mathbb{Z}},s)$ be the compact topological dynamic consists of the space of all binary $\mathbb{Z}$-sequences equipped with product topology on which operator shift $s$ is acting. This system is called  
Bernoulli system. 
For every finite binary sequence $I$, by $[I]$ we mean the set of those $J \in 2^{\mathbb{Z}}$ containing $I$ in its initial segment, more precisely $J(i)=I(i)$ for every $i=1,\ldots,|I|$. We call $[I]$ a basic open set in $2^{\mathbb{Z}}$. The family of basic open sets for a basis for the topology. Note that topology of this system is metrizable and coming from the metric $d((a_i)_{i \in \mathbb{Z}},(b_i)_{i \in \mathbb{Z}})=\frac{1}{2^t}$ where $t=\min\{|i|: a_i \not = b_i, i \in \mathbb{Z}\}$.
Also note that $d$ is an ultrametric (i.e for every $x,y,z$ we have that $d(x,y) \leqslant \max \{d(x,z),d(y,z)\}$). It is also a non Archimedean metric space and for each ball every point is a center.  
Bernoulli system is an important example in symbolic dynamical systems and its dynamic is chaotic.
For every $0<p<1$, by a $p$-Bernoulli measure on $\mathcal{B}$ we mean the unique measure denoted by $\mu$ such that for every finite binary sequence $I$, we have $\mu([I])=p^{r_I}(1-p)^{S_I}$ where $r_I$ and $s_I$ denote the number of $1$'s and $0$'s appearing in $I$ respectively. One can see that $p$-Bernoulli measures are shift invariant on $\mathcal{B}$.

\vspace{1.5mm}

\noindent \textbf{{\large Some notations}}

\vspace{1.5mm}

We set some notations. Let $\phi$ be a formula.
We use $\phi^0$ and $\phi^1$ for the negation of $\phi$ and $\phi$ respectively.
Also let $\bar{\phi}(x,y):=\phi(y,x)$.
We denote by $\bar{1}\bar{0}$ the $\mathbb{Z}$-sequence $I$ with $I(i)=1$ for every $i \leqslant 0$ and $I(i)=0$ for every $i >0$.
Similarly, $\bar{0}\bar{1}$ is defined in the converse way.
In fact we can identify $\bar{1}\bar{0}$ and $\bar{0}\bar{1}$ with non positive and non negative parts of $\mathbb{Z}$ respectively.
By ${(\bar{1}\bar{0})}_n$ we mean finite part of $\bar{0}\bar{1}$ from $-n$'th to $n$'th coordinates. ${(\bar{0}\bar{1})}_n$ is also defined similarly.
We say that a binary sequence is universal if it contains every finite binary sequence as a subsequence. Note that in this paper we frequently identify subsets of $\mathbb{N}$ and $\mathbb{Z}$ with binary $\mathbb{N}$-sequences and $\mathbb{Z}$-sequences respectively.

Assume that $I$ is an infinite binary sequence and $a \in \mathbb{Z}$. By $I+a$ we mean the sequence obtained from $I$ after $a$ shifts. If $I$ is a finite sequence then $I'=I+a$ is same as $I$ but starts from $a+1$ (more precisely it can be seen as $I':[a+1,\ldots,a+|I|] \rightarrow \{0,1\}$ with $I'(i)=I(i-a)$ for every $i \in [a+1,\ldots,a+|I|]$).
Also for any binary sequence $I$, by $I^c$ we mean the sequence obtained by replacing $1$ and $0$ in $I$ to each other. Let $\{I_{\alpha}: \alpha \in \Omega \}$ be a family of binary $\mathbb{Z}$-sequences. By $\bigcap_{\alpha \in \Omega} I_{\alpha}$ we mean a binary $\mathbb{Z}$-sequence $J$ where for every $i \in \mathbb{Z}$ we have $J(i)=1$ if and only if $J_{\alpha}(i)=1$ for every $\alpha \in \Omega$.

Let $M$ be a model and $a \in M^n$ for some $n$. We define $W_{a}(A)$ to be the set of auotomorphisms $\sigma$ of $M$ with the property that $\sigma$-orbit of $a$ is an $A$-indiscernible sequence. When $A=\emptyset$, we denote it by $W_{a}$.
For an automorphism $\sigma$ by $\sigma$-orbit of an element $\a$ we mean the set $\mathcal{O}_{\sigma}(a)=\{\sigma^n(a): n \in \mathbb{Z}\}$.
\vspace{2mm}

\noindent \textbf{{\large Some operations on $\mathcal{B}$}}

\vspace{2mm}

Let $I, J$ be two elements of $\mathcal{B}$. We say that $I$ is a switching of $J$ if $I$ is obtained from $J$ by replacing some $01$ to $10$ or replacing some $10$ to $01$ in $J$. Obviously if $I$ be a switching of $J$ then $J$ is a switching of $I$. 
One can consider $\mathcal{B}$ as a graph where two point are connected if one is switching of the other one. 
For a subset $H$ of $\mathcal{B}$, we denote by $\langle H \rangle_0$ the set containing $H$ and all $I \in \mathcal{B}$ so that there exists some finite path between some element of $H$ and $I$.
In fact looking to $\mathcal{B}$ in this way, $\langle H \rangle_0$ would be the union of all connected components of $\mathcal{B}$ which have nonempty intersection with $H$.
Let $w:\mathcal{B} \rightarrow \mathcal{B}$ be the function that for every $I$, switches the values of $I(0)$ and $I(1)$ to each other. Obviously $w$ is a continuous function.
We call a subset of $\mathcal{B}$ SW-closed if it is closed under the action of shift and the function $w$. Also we define the SW-closure of a set $H \subseteq \mathcal{B}$ to be the minimal SW-closed set containing $H$. We denote the SW-closure of a set $H$ by $\langle H \rangle$. One can see that $\langle H \rangle$ is the closure of $\langle H \rangle_0$ under the shift action.
Note that SW-closeness implies closeness under any switching. This is because switching in the $k$'th coordinate can be written as $s^{k}(w(s^{-k}(I)))$ where $s$ is the shift function.

\begin{rem} \label{<01>dense}
\emph{
$\langle \{\bar{1}\bar{0}\} \rangle$ and $\langle \{\bar{0}\bar{1}\} \rangle$ are dense in $\mathcal{B}$ with respect to the natural metric on $\mathcal{B}$.
}
\end{rem}

\begin{dfn}\label{defsymbrepandtotalsymbmap}
Let $M$ be a model.
\emph{\begin{enumerate}
\item{Let $\sigma \in Aut(M)$ and $U \subseteq M$. 
We recall the Definition \ref{defsymbimagsingleelement} in this case and define the \textit{symbolic representation} with respect to $\sigma$ and $U$ as the following map
$$\xi_{\sigma,U}: M \rightarrow 2^{\mathbb{Z}}$$
$$a \rightarrow \{n \in \mathbb{Z} : \sigma^n(a) \in U\}.$$
}
\item{
Let $\phi(x,y)$ be a formula, $A \subseteq M$, $G=Aut(M/A)$ and $\mathcal{U}=\{\phi(x,b): b \in M\}$. 
We define the \textit{$\mathbb{Z}$-total symbolic image} of the formula $\phi(x,y)$ on parameter set $A$ with respect to the action of $G$ to be 
$\rho_{G,\mathbb{Z},\mathcal{U}}(M):=\bigcup_{U \in \mathcal{U}} \rho_{G,\mathbb{Z},U}(M)$
where 
$\rho_{G,\mathbb{Z},U}(M):=\bigcup_{\sigma \in G} \xi_{\sigma,U}(M)$. 
When $A=\emptyset$ we use notation $\rho_{\phi}(M)$ and call it the 
\textit{$\mathbb{Z}$-total symbolic image} of the formula $\phi(x,y)$.}\label{deftotalsymbimg}
\end{enumerate}}
\end{dfn}

The following proposition roughly says that for a definable set no essentially new element in symbolic images
will be added by going to elementary extensions.

\begin{pro}\label{symbimginelemext}
Let $M$ be a model, $b \in M$ and $U=\phi(M,b)$ be a definable set in $M$. Let $\sigma$ be a map 
from $M$ to $M$ and $N$ be an elementary extension of $M$.
Also let $\tau$ be an extension of $\sigma$ to $N$ and $U^N=\phi(N,b)$. Then
\begin{enumerate}
\item{
The set $\xi_{\sigma,U}(M)$ is dense in $\xi_{\tau,U^N}(N)$.
}\label{symbolicimagesmallmodelisdenseinbigmodel}
\item{
If $N$ is $\aleph_1$-saturated and $\sigma$ and $\tau$ are automorphisms 
of $M$ and $N$ respectively, then $\xi_{\tau,U^N}(N)$ is equal to $cl(\xi_{\sigma,U}(M))$, the topological closure of $\xi_{\sigma,U}(M)$ in $\mathcal{B}$.
}\label{symbolicimagebigmodelisistopclosureofsmallmodel}
\end{enumerate}
\end{pro}
\proof
\ref{symbolicimagesmallmodelisdenseinbigmodel})
Assume for contradiction that there exists $I \in \xi_{\tau,U^N}(N) \setminus cl(\xi_{\sigma,U}(M))$. Since $I \not \in cl(\xi_{\sigma,U}(M))$, for some initial segment of $I$, say $J=(I(-n),\ldots,I(n))$, no element of $\xi_{\sigma,U}(M)$ contains $J$ as initial segment. 
Let $b_i= \tau^i(b)$ for each $i \in \mathbb{Z}$. Since $\tau$ is an extension of $\sigma$ and $b \in M$, then $b_i$'s are in $M$. So  
$M \models *$ where $*$ is the following sentence
$$*: \not \exists x (\phi^{I(-n)}(x,b_n) \wedge \phi^{I(-n+1)}(x,b_{n-1}) \wedge \ldots \wedge \phi^{I(n)}(x,b_{-n})).$$
Since $M \preceq N$ and $*$ is a sentence in the language $\mathcal{L}(M)$, we have 
$N \models *$. This implies that there is no element in $\xi_{\tau,U^N}(N)$ containing $J$ as initial segment. This contradicts with the existence of $I$. So $\xi_{\tau,U^N}(N) \subseteq cl(\xi_{\sigma,U}(M))$.

\vspace{2mm}

\ref{symbolicimagebigmodelisistopclosureofsmallmodel})
Using part \ref{symbolicimagesmallmodelisdenseinbigmodel}, it is enough to show that $cl(\xi_{\sigma,U}(M)) \subseteq \xi_{\tau,U^N}(N)$.
For every natural number $j$ let $d_j$ be an element of $M$ and $I_j:=\xi_{\sigma,U}(d_j)$.
Assume that $I \in \mathcal{B}$ be such that $I_j$'s converge to $I$. 
Let $c_i=\sigma^i(b)$ for each $i \in \mathbb{Z}$. 
We define a partial type over 
$\mathcal{O}_{\sigma}(b)$ as follows. 
$$p=\{\phi^{I(i)}(x,c_{-i}), i \in \mathbb{Z}\}.$$
We recall that by $\phi^0$ and $\phi^1$ we mean $\neg \phi$ and $\phi$ respectively.
One can see that $p$ is consistent with $T$ and is a partial type. 
Now by $\aleph_1$-saturation of $N$ and since the parameters of $p$ are in $M$, $p$ is realized in $N$ by some $d \in M$.
Since $\tau$ is an automorphism, for every arbitrary fix $i \in \mathbb{Z}$ we have that
$$\tau^i(d) \in U^N \Leftrightarrow d \in \tau^{-i}(\phi(M,b)) \Leftrightarrow \phi(d,\tau^{-i}(b)) \Leftrightarrow \phi(d,\sigma^{-i}(b)) \Leftrightarrow \phi(d,c_{-i}) \Leftrightarrow I(i)=1$$
So we have that $\xi_{\tau,U^N}(d)=I$.
Thus $I \in \xi_{\tau,U^N}(N)$.
So $cl(\xi_{\sigma,U}(M)) \subseteq \xi_{\tau,U^N}(N)$ and the proof is complete.
$\ \ \square$
\begin{cor}\label{aleph1satimpliesCloseness}
Let $M$ be a $\aleph_1$-saturated model of a theory $T$ and $U=\phi(M,b)$. For every $\sigma \in Aut(M)$ 
the followings hold.
\begin{enumerate}
\item{$\xi_{\sigma,U}(M)$ is closed in $\mathcal{B}$.}\label{saturationimpliesclose}

\item{For every $M \preceq N$ and every $\tau \in Aut(N)$
extending $\sigma$, we have that $\xi_{\tau,V}(N)=\xi_{\sigma,U}(M)$ where $V=\phi(N,b)$.}\label{saturatedimpliessameinextensions}
\end{enumerate}
\end{cor}
\proof
\ref{saturationimpliesclose})
Using part \ref{symbolicimagebigmodelisistopclosureofsmallmodel} of \ref{symbimginelemext} and letting $M=N$.

\ref{saturatedimpliessameinextensions})
Using parts 
\ref{symbolicimagesmallmodelisdenseinbigmodel} and \ref{symbolicimagebigmodelisistopclosureofsmallmodel} of \ref{symbimginelemext}, $\xi_{\sigma,U}(M)$ is dense in $\xi_{\tau,V}(N)$ and also is closed. Therefore $\xi_{\tau,V}(N)=\xi_{\sigma,U}(M)$. $\ \ \square$

\begin{rem}
\emph{If $M$ is not  $\aleph_1$-saturated then $\xi_{\sigma, U}(M)$ might not be closed in $\mathcal{B}$. For example let $M$ be the countable random graph (the Rado graph), $a \in M$, $U:=\{x \in M: x \thicksim a\}$ be a definable set (where $\thicksim$ denotes the adjacency relation in graph) and $\sigma$ be a cyclic automorphism. Note that in this case one can look at $M$ as the Cayley graph of $(\mathbb{Z},+,S)$ for a suitable universal sequence $S$ and automorphism $\sigma$ as the shift on $\mathbb{Z}$.
Using this point of view $\xi_{\sigma, U}(M)$ consists of the shift orbit of the sequence $\xi_{\sigma, U}(a)$. So $\bar{1} \not \in \xi_{\sigma, U}(M)$. But since $S$ is a universal $\mathbb{Z}$-sequence, $\xi_{\sigma, U}(M)$ is dense in $2^{\mathbb{Z}}$. So it is not closed.}
\end{rem}

\noindent  \textbf{Push forward measures via symbolic representations}

\vspace{1mm}

Let $M$ be a model, $U$ be definable set of some instance of some formula with parameter from $M$, $\mu$ be a Keisler measure on $S(M)$ and $\sigma$ be an automorphism.
Obviously $\sigma$ acts on $S(M)$ and $U$ can be seen as a subset of $S(M)$. One may consider symbolic representation with respect to $\sigma$ and $U$ in $S(M)$.
The following statement shows that symbolic representations are continuous function on spaces of types.

\begin{lem}\label{continuousandcomutativityrhosigmashift}
\begin{enumerate}
\item{$\xi_{\sigma, U}$ is a continuous function.}\label{symbrepiscont}
\item{$\xi_{\sigma, U} \ o \ \sigma = S^{-1} \ o \ \xi_{\sigma, U}$ and 
$\xi_{\sigma, U} \ o \ \sigma^{-1} = S \ o \ \xi_{\sigma, U}$}\label{symbrepcommute}
\end{enumerate}
\end{lem}
\proof
$\ref{symbrepiscont})$ For definable set $U$, by $U^1$ and $U^0$ we mean $U$ and $U^c$ respectively. Let $A=\xi_{\sigma, U}(S(M))$. Also let $I$ be an arbitrary element of $2^{\mathbb{Z}}$. We show that for every $r \in \mathbb{R}^+$, $\xi_{\sigma, U}^{-1}(B_r(I))$ is open in $S(M)$ where $B_r(I)$ is the open ball with radius $r$ and center $I$ in $\mathcal{B}$. There exists some $n_r \in \mathbb{N}$ such that $B_r(I)$
consists of those $J \in 2^{\mathbb{Z}}$ such that $I(i)=J(i)$ for every $-n_r \leqslant i \leqslant n_r$.
So we have that
$$\xi_{\sigma, U}^{-1}(B_r(I))=\{p \in S(M): U^{I(i)} \in \sigma^i(p) \ for \ every -n_r \leqslant i\leqslant n_r \}
=\bigcap_{i=-n_r}^{n_r} \sigma^i( U^{I(i)}).$$
But the last expression is an open set in $S(M)$. So we are done.

\vspace{2mm}

$\ref{symbrepcommute})$
For simplicity we use the notation $\xi$ for $\xi_{\sigma, U}$.
Let $p \in S(M)$. So $$\xi \ o \ \sigma(p)=\{i: \sigma^i(\sigma(p)) \in U\}=\{i-1: \sigma^i(p) \in U\}=\xi(p)-1=s^{-1} \ o \ \xi(p).$$
Similarly $$\xi \ o \ \sigma^{-1}(p)=\{i: \sigma^i(\sigma^{-1}(p)) \in U\}=\{i+1: \sigma^i(p) \in U\}=\xi(p)+1=s \ o \ \xi(p). \ \ \square$$

\vspace{1mm}

Let $\mu$ be a measure on $S(M)$. By previous lemma $\xi_{\sigma,U}$ is a continuous function on $2^{\mathbb{Z}}$. Hence it is a measurable function and can push forward the measure $\mu$ from $S(M)$ on $2^{\mathbb{Z}}$. We denote this pushed forward measure by $\nu_{\mu,\sigma,U}$. Note that $\nu_{\mu,\sigma,U}$ concentrates on $\xi_{\sigma,U}(S(M))$. 
The following lemma shows that if moreover $\mu$ is $\sigma$-invariant, then 
$\nu_{\mu,\sigma,U}$ is $s^{-1}$-invariant where $s$ and $s^{-1}$ are right and left shifts on $2^{\mathbb{Z}}$ respectively. Moreover if $\mu$ is both $\sigma$ and $\sigma^{-1}$ invariant then induced measure will be two sided shift-invariant. 
Note that shift invariant measures on $2^{\mathbb{Z}}$ are extensively studied in the literature.

\begin{rem}\label{ASRpushforwardmeasureisinvariant}
\emph{
If $\mu$ on $S(M)$ is $\sigma$-invariant ($\sigma^{-1}$-invariant) then $\nu_{\mu,\sigma,U}$ on $2^{\mathbb{Z}}$ will be $s^{-1}$-invariant ($s$-invariant).
}
\end{rem}
\proof
We use $\xi$ instead of $\xi_{G,U}$ and $\nu$ instead of $\nu_{\mu,\sigma,U}$.
Obviously, for every measurable $A \subseteq 2^{\mathbb{Z}}$, we have that $\nu(A)=\nu(\xi \ o \ \xi^{-1}(A))$.
Using part \ref{symbrepcommute} of Remark \ref{continuousandcomutativityrhosigmashift}
we have that $\xi^{-1}(s^{-1}(A))=\sigma (\xi^{-1}(A))$ and 
$\xi^{-1}(s(A))=\sigma^{-1} (\xi^{-1}(A))$.
Now if $\mu$ is $\sigma^{-1}$-invariant one has  
$$\nu(s^{-1}(A))=\mu(\xi^{-1}(s^{-1}(A))) =\mu((\sigma (\xi^{-1}(A))))=\mu((\xi^{-1}(A)))=\nu(A).$$
Similarly if $\mu$ is $\sigma$-invariant one has that
$$\nu(s(A))=\mu(\xi^{-1}(s(A))) =\mu((\sigma^{-1} (\xi^{-1}(A))))=\mu((\xi^{-1}(A)))=\nu(A). \ \ \square$$

\subsubsection{Independence property}

The following statement characterizes independence property in terms of symbolic images of formulas.
\begin{thm} \label{symboliccharacterizationofNIP}  
Let $T$ be a theory and $\phi(x,y)$ a formula. Then the followings are equivalent.
\begin{enumerate}
\item{The formula $\phi(x,y)$ has IP.}\label{phihasIP}
\item{ 
There exists some model $M$ of $T$, some instance $U=\phi(M,b)$ for some $b \in M$ and some automorphism $\sigma$ of $M$ such that $\xi_{\sigma,U}(M)$ is dense in $\mathcal{B}$.}\label{imageisDense}
\item{ 
There exists some model $M$ of $T$ and some automorphism $\sigma$ of $M$ such that for every $n \in \mathbb{N}$,
there exists some instance $U_n=\phi(M,b_n)$ for some $b_n \in M$ such that every binary sequence of lenght $n$ appears in some element of $\xi_{\sigma,U}(M)$ as a subsequence.}\label{imagehasfiniteuniversalseq}
\item{ 
There exists some model $M$ of $T$, some instance $U=\phi(M,b)$ for some $b \in M$ and some automorphism $\sigma$ of $M$ such that $\xi_{\sigma,U}(M)$ contains a universal sequence.}\label{imagehasUniversal}
\item{ 
There exists some model $M$ of $T$, some instance $U=\phi(M,b)$ for some $b \in M$ and some automorphism $\sigma$ of $M$ such that $\xi_{\sigma,U}(M)=\mathcal{B}$.}\label{mapissurjective}
\item{There exists some model $M$ of $T$ such that $\rho_{\phi}(M)$ is dense.}
\label{totalsymbimgdense}

\end{enumerate}
\end{thm}
\proof
($\ref{imageisDense} \Rightarrow \ref{phihasIP}$)
Let $n \in \mathbb{N}$ be arbitrary. We show that there exists a witness for IP of length $n$.
Let $J=J_1,\ldots,J_{2^n}$ be the sequence of length $n.2^n$ obtained by concatenation of all
binary sequences $J_i$'s of length $n$.
Since $\xi_{\sigma, U}(M)$ is dense in $\mathcal{B}$, there exists some $x_J \in M$ such that 
$\{i \ : \ 1 \leqslant i \leqslant 2^n, \ \sigma^i(x_J) \in U \}$ is same as $J$ when we regard it as a binary sequence.
We claim that the sequence $(\sigma^i(x_J))_{i=1,\ldots,n}$ 
witnesses IP of length $n$.
Let $\{I_1,I_2\}$ be an arbitrary partitioning of $I=\{1,\ldots,n\}$.
There is exactly one $1 \leqslant t \leqslant 2^n$ such that $J_t(j)=1$ for $j \in I_1$ and $J_t(j)=0$ for $j \in I_2$.
So $J_t$ corresponds to $n(t-1)+1$'th up to $nt$'th coordinates of $J$.
So for each $i \in I_1+n(t-1)$ we have that
$\phi(\sigma^i(x_J),b)$. Also for each $i\in I_2+n(t-1)$ we have that $\neg \phi(\sigma^i(x_J),b)$.
Let $d=\sigma^{-n(t-1)}(b)$.
Therefore we have 
$\phi(\sigma^i(x_J),d)$ for every $i \in I_1$ and
$\neg \phi(\sigma^i(x_J),d)$ for every $i \in I_2$.
So we have witnessed IP of length $n$.

\vspace{1.5mm}

$(\ref{imagehasUniversal} \Rightarrow \ref{imageisDense})$
Let $x_I \in M$ be such that $I:=\xi_{\sigma, U}(x_I)$ be a universal sequence.
We denote the sequence $I+n$ (sequence obtained from $I$ after $n$ shifts) by $I_n$.
So for every $n \in \mathbb{Z}$, we have $\xi_{\sigma, U}(\sigma^n(x_I))=I_n$. 
Hence $I_n \in \xi_{\sigma, U}(M)$.
So for every finite sequence of length $2n$ ($n \in \mathbb{N}$), there exists some element of $\xi_{\sigma, U}(M)$ containing that sequence in $-n$'th up to $n$'th coordinates. This implies that $\xi_{\sigma, U}(M)$ is dense in $\mathcal{B}$.

($\ref{imagehasfiniteuniversalseq} \Rightarrow \ref{phihasIP}$)
For every $n \in \mathbb{N}$, using an argument similar to that of part 
$(\ref{imageisDense} \Rightarrow \ref{phihasIP})$, IP is witnessed with length $n$.

($\ref{totalsymbimgdense}  \Rightarrow \ref{phihasIP}$)
Let $n \in \mathbb{N}$ be arbitrary and $J$ be as in the argument of part 
$(\ref{imageisDense} \Rightarrow \ref{phihasIP})$.
By part \ref{deftotalsymbimg} of the Definition \ref{defsymbrepandtotalsymbmap}, we have 
$\rho_{\phi}(M)=\bigcup_{b \in M} \bigcup_{\sigma \in Aut(M)} \xi_{\sigma,U_b}(M)$ where 
$U_b=\phi(x,b)$ for every $b$.
Since $\rho_{\phi}(M)$ is dense there exists some $\sigma \in Aut(M)$, some $b \in M$ and some $x_J \in M$ such that $\{i : \ 0 < i \leqslant n, \ \sigma^i(x_J) \in U_b \}=J$.
Now again similar to the argument of part $(\ref{imageisDense} \Rightarrow \ref{phihasIP})$, 
$\{\sigma^1(x_J), \ldots, \sigma^n(x_J)\}$ witnesses IP with length $n$ for $\phi$.

\vspace{1.5mm}

$(\ref{mapissurjective} \Rightarrow \ref{imageisDense})$,
$(\ref{mapissurjective} \Rightarrow \ref{imagehasUniversal})$, $(\ref{imagehasUniversal} \Rightarrow \ref{imagehasfiniteuniversalseq})$ and 
$(\ref{imageisDense} \Rightarrow \ref{totalsymbimgdense})$
are obvious.

\vspace{1.5mm}

$(\ref{phihasIP} \Rightarrow \ref{mapissurjective})$
By a theorem of Shelah (see for example Theorem 12.18 of \cite{Poizatmodeltheory}) 
there exists a formula 
$\phi(x,y)$ with $|y|=1$ witnessing IP. Also by some fact (See for example Lemma 12.16 of \cite{Poizatmodeltheory}), the formula $\psi$ defined with $\psi(x,y):=\phi(y,x)$ has IP too.
By saturation, there exists an indiscernible sequence $\{a_i\}_{i \in \mathbb{Z}}$ witnessing IP for $\phi$. Again by saturation there exists some automorphism of $M$, say $\sigma$, such
$a_i=\sigma^i(a_0)$. Let $I$ be an arbitrary binary $\mathbb{Z}$-sequence.
By IP, there exists $b \in M$ such that for every $i \in \mathbb{Z}$
we have $\phi(a_i,b)$ if and only if $i \in -I$.
Now let
$b_i=\sigma^i(b)$.
Then clearly we have $\psi(b_i,a_0)$ if and only if $\phi(a_{-i},b)$.
By letting $U=\psi(M,a_0)$, we have that $\xi_{\sigma,U}(b)=I$. Now the proof is complete.
$\square$

\vspace{1.5mm}

One sees that the equivalence of parts $\ref{phihasIP}$ and $\ref{totalsymbimgdense}$ in the above theorem shows that whether or not $\phi$ has IP, is reflected in the 
$\rho_{\phi}(M)$, the $\mathbb{Z}$-total symbolic image of $\phi$.

Now we give a characterization of NIP property in terms of Bernoulli measures and push-forwarded measure on $\mathcal{B}$.

\begin{pro}\label{IPiffBernoullimeasonsymbimage} 
Let $T$ be a theory and $\phi(x,y)$ a formula. Then the followings are equivalent.
\begin{enumerate}
\item{The formula $\phi(x,y)$ has IP.}\label{IPiffBernoul-phihasIP}
\item{There exists some model $M$ of $T$, some instance $U=\phi(M,b)$ for some $b \in M$, some automorphism $\sigma$ of $M$ and some $\sigma$-invariant measure 
$\mu$ on $S(M)$ (or $M$) such that the push forward measure $\mu_{\xi}=\xi_{\sigma,U}(\mu)$ is the $p$-Bernoulli measure 
on $\mathcal{B}$ for some $0<p<1$.}\label{IPiffBernoul-imagemeasissomeBernoulli}
\item{ 
For any given $0<p<1$, there exists some model $M$ of $T$, some instance $U=\phi(M,b)$ for some $b \in M$, some automorphism $\sigma$ of $M$ and some $\sigma$-invariant measure 
$\mu$ on $S(M)$ (or $M$) such that the push forward measure $\mu_{\xi}=\xi_{\sigma,U}(\mu)$ is the $p$-Bernoulli measure on $\mathcal{B}$.}\label{IPiffBernoul-imagemeasisgivenBernoulli}
\end{enumerate}
\end{pro}
\proof
$(\ref{IPiffBernoul-phihasIP} \Rightarrow \ref{IPiffBernoul-imagemeasisgivenBernoulli})$
Assume that $0<p<1$ is given and $\phi(x,y)$ has IP. Let $M$ be a saturated enough model of $T$ and let $B:=(a_i)_{i \in \mathbb{Z}}$ be an indiscernible sequence witnessing IP for $\phi$. 
Because of saturation, there exists some $\sigma \in Aut(M)$ such that $\sigma^i(a_0)=a_i$ for every $i$. By Remark \ref{existenceofpnormalnumbers}, there exists some $p$-normal number $u$.
Let $J$ be 
the $\mathbb{Z}$-sequence associated to $u$
(defined in the Definition \ref{sequencefromrealnmber}).
Since $B$ witnesses IP, there exists some $b \in M$ such that 
$\phi(a_i,b) \Leftrightarrow i \in J$. Let $U=\phi(M,b)$.
So $a_i=\sigma^i(a_0) \in U^{J(i)}$ for every $1 \leqslant i \leqslant |J|$.
Let $\mathcal{F}$ be an ultrafilter 
extending the Frechet filter and let $\mu^{\sigma,a_0}_{\mathcal{F}}$ be the limit frequency measure on $Def(M)$
with respect to $B$ defined in the Definition \ref{dfnlimitfrequencymeasures}. So
for every definable set $D$, 
we have $\mu^{\sigma,a_0}_{\mathcal{F}}(D)=dns_{\mathcal{F}}(\xi_{\sigma,D}(a_0))$.
As explained in the Remark \ref{sigmaaddivitemeasureonmodel}, this measure could be seen as a measure on $S(M)$.
Now we claim that the pushed forward measure $\nu$ induced by $\mu^{\sigma,a_0}_{\mathcal{F}}$ on $2^{\mathbb{Z}}$ via the map $\xi_{\sigma,U}$ is a $p$-Bernoulli measure. 
For showing that, it is enough to show that the measures of cylindrical sets are as what the $p$-Bernoulli measure gives to these sets. Let $W$ be an arbitrary finite binary sequence
and $\alpha$ and $\beta$ be number of appearance of $1$'s and $0$'s in $W$ respectively.
By using of notations we have defined earlier we have that
$$\nu([W])=
\mu^{\sigma,a_0}_{\mathcal{F}}(\{x \in M: \xi_{\sigma,U}(x) \in [W]\})=
\mu^{\sigma,a_0}_{\mathcal{F}}(\{x \in M: \sigma^j(x) \in U^{J(j)}, 1 \leqslant j \leqslant |W|\})$$
$$=\mu^{\sigma,a_0}_{\mathcal{F}}(\bigcap_{1 \leqslant j \leqslant |W|} \sigma^{-j}(U^{W(j)}))
=dns_{\mathcal{F}}(\{i \in \mathbb{Z}: \sigma^i(a_0) \in \bigcap_{1 \leqslant j \leqslant |W|} \sigma^{-j}(U^{W(j)})\})$$
$$
=dns_{\mathcal{F}}(\bigcap_{1 \leqslant j \leqslant |W|}\{i \in \mathbb{Z}: \sigma^i(a_0) \in \sigma^{-j}(U^{W(j)})\})
=dns_{\mathcal{F}}(\bigcap_{1 \leqslant j \leqslant |W|}\{i \in \mathbb{Z}: \sigma^{i+j}(a_0) \in U^{W(j)}\})$$
$$=dns_{\mathcal{F}}(\bigcap_{1 \leqslant j \leqslant |W|}(\xi_{\sigma,U^{W(j)}}(a_0)-j))
=dns_{\mathcal{F}}(\bigcap_{1 \leqslant j \leqslant |W|} ( (J-j)^{W(i)}).$$
Now using Remark \ref{capI-j=<W,I>} we have that 
$$\nu([W])=dns_{\mathcal{F}}(\langle W,J \rangle).$$
But since $J$ is the sequence associated to a $p$-normal number we have that
$$dns_{\mathcal{F}}(\langle W,J \rangle)=p^{\alpha}(1-p)^{\beta}.$$
So claim is proved and \ref{IPiffBernoul-imagemeasisgivenBernoulli} holds.

\vspace{1.5mm}

$(\ref{IPiffBernoul-imagemeasisgivenBernoulli} \Rightarrow \ref{IPiffBernoul-imagemeasissomeBernoulli})$ is obvious.

\vspace{1.5mm}

$(\ref{IPiffBernoul-imagemeasissomeBernoulli} \Rightarrow \ref{IPiffBernoul-phihasIP})$
Assume that \ref{IPiffBernoul-imagemeasissomeBernoulli} holds. So there exists some $0<p<1$ and some $p$-Bernoulli measure $\mu_{\xi}$ which is  
pushed forward from some measure $\mu$ on $M$.
For every $i \not = 0$ we have that 
$$\mu(U \cap \sigma^i(U))=\mu_{\xi}(\{I: I(1)=1,I(i)=1\})=p^2.$$
So
$$\mu(U \triangle \sigma^i(U))=\mu(U)+\mu(\sigma^i(U))-2\mu(U \cap \sigma^i(U))
\geqslant 2(p-p^2).$$
Now by using the Theorem \ref{K-NIPKeislermeasures} the formula $\phi(\x,\y)$ has IP.
$\square$

\subsubsection{Order property}
The following statement characterizes order property in terms of symbolic images of formulas.
\begin{pro} \label{symboliccharacterizationofOP}
Let $T$ be a theory and $\phi(x,y)$ be a formula. Then the followings are equivalent.

\begin{enumerate}
\item{
The formula $\phi(x,y)$ has order property.} \label{phihasOP} 

\item{
There exists some model $M$ of $T$, some instance $U=\phi(M,b)$ for some $b \in M$, some automorphism $\sigma$ of $M$ and some $a \in M$ such that $\xi_{\sigma,U}(a)=\bar{1}\bar{0}$.} \label{10imageofsomea}

\item
{Same as \ref{10imageofsomea} with additional property that $\sigma \in W_{b}$.} \label{10imageofsomeawithtauinWb}

\item{
There exists some model $M$ of $T$, some instance $U=\phi(M,b)$ for some $b \in M$, some automorphism $\sigma$ of $M$ such that for every $n \in \mathbb{N}$, there exists some element in $\xi_{\sigma,U}(M)$ which contains $({\bar{1}\bar{0}})_n$  
as a subsequence.} \label{long10inimage}

\item{
There exists some model $M$ of $T$ and some automorphism $\sigma$ of $M$ such that for every 
$n \in \mathbb{N}$,
there exists some instance $U_n=\phi(M,b_n)$ for some $b_n \in M$, 
such that there exists some element in $\xi_{\sigma,U_n}(M)$ which contains
$({\bar{1}\bar{0}})_n$ as a subsequence.} \label{everylong10inimageofsome}

\item{
There exists some model $M$ of $T$, some instance $U=\phi(M,b)$ for some $b \in M$, some automorphism $\sigma$ of $M$ such that the topological closure of $\xi_{\sigma, U}(M)$ contains  $\bar{1}\bar{0}$. 
}\label{orbitclosurecontains10}

\item{There exists some model $M$ of $T$ such that $\bar{1}\bar{0} \in \rho_{\phi}(M)$.}\label{totalsymbimgcontains10}

\end{enumerate}
\end{pro}
\proof
$(\ref{phihasOP} \Rightarrow \ref{10imageofsomea})$ Using the Lemma \ref{differentequivalencesofOP}, there exists some model $M$, some indiscernible sequence $I=(a_i)_{i \in \mathbb{Z}}$ and a sequence $J=(b_j)_{i \in \mathbb{Z}}$ in $M$ such that $\phi(a_i,b_j) \Leftrightarrow i \leqslant j$ for every $i,j \in \mathbb{Z}$.
We may assume that $M$ is saturated enough. By indiscerniblity of $I$, there exists some automorphism $\sigma$ of $M$ such that $\sigma^i(a_0)=a_i$ for every $i \in \mathbb{Z}$. So letting $U=\phi(M,b_0)$ we have that 
$$\xi_{\sigma,U}(a_0)=\{i \in \mathbb{Z}, \phi(a_i,b_0)\}=\bar{1}\bar{0}.$$

$(\ref{long10inimage} \Rightarrow \ref{phihasOP})$ 
For every $n \in \mathbb{N}$, we find a witness for order property with length $n$.
By assumption and also the fact that $\xi_{\sigma,U}(M)$ is shift invariant, there exists some $a \in M$ such that $\xi_{\sigma,U}(a)$ contains ${(\bar{1}\bar{0})}_{n}$ in the coordinates $-n$ up to $n$. So we have that
$\{-n \leqslant i \leqslant n: \sigma^i(a) \in \phi(M,b)\}={(\bar{1}\bar{0})}_{n}$. 
Let $a_i=\sigma^i(a)$ and $b_i=\sigma^i(b)$ for every $i \in \mathbb{Z}$. Then for every $0 \leqslant i,j \leqslant n$ we have that
$\phi(a_i,b_j) \Leftrightarrow \phi(a_{i-j},b) \Leftrightarrow i \leqslant j$.  
So $\phi$ has OP of length $n$.

\vspace{1.5mm}

$(\ref{totalsymbimgcontains10} \Rightarrow \ref{phihasOP})$
By part \ref{deftotalsymbimg} of the Definition \ref{defsymbrepandtotalsymbmap} we have 
$\rho_{\phi}(M)=\bigcup_{\b \in M} \bigcup_{\sigma \in Aut(M)} \xi_{\sigma,U_{b}}(M)$
where $U_{b}=\phi(x,b)$ for every $b$.
So there exists some $\sigma \in Aut(M)$, some $b \in M$ and some $a \in M$ such that 
$\{i \in \mathbb{Z}, \ \sigma^i(a) \in \phi(M,b)\}={\bar{1}\bar{0}}$.
Let $a_i=\sigma^i(a)$ and $b_i=\sigma^i(b)$ for every $i \in \mathbb{Z}$. So we have that $\phi(a_i,b_j) \Leftrightarrow \phi(a_{i-j},b)$ if and only if $i \leqslant j$
for every $i,j \in \mathbb{Z}$.

\vspace{1.5mm}

$(\ref{10imageofsomea} \Rightarrow \ref{10imageofsomeawithtauinWb})$ By Ramsey theorem and compactness.
Also $(\ref{everylong10inimageofsome} \Rightarrow \ref{orbitclosurecontains10})$ is by compactness. 

$(\ref{10imageofsomea} \Rightarrow \ref{long10inimage})$ and $(\ref{10imageofsomea} \Rightarrow \ref{totalsymbimgcontains10})$ are obvious.

$(\ref{long10inimage} \Rightarrow \ref{orbitclosurecontains10})$, $(\ref{orbitclosurecontains10} \Rightarrow \ref{long10inimage})$, $(\ref{10imageofsomeawithtauinWb} \Rightarrow \ref{10imageofsomea})$ and $(\ref{orbitclosurecontains10} \Rightarrow \ref{everylong10inimageofsome})$ are easy.
$\ \ \square$

\vspace{2mm}

Note that assuming $M$ is $\aleph_0$-saturated, then by Proposition \ref{symboliccharacterizationofOP} and the Corollary \ref{aleph1satimpliesCloseness},
$\phi$ has OP if and only if for some instance $U$ of $\phi$ and some automorphism $\sigma$ of $M$ we have that $\bar{1}\bar{0} \in \xi_{\sigma, U}(M)$. One sees that the equivalence of $\ref{phihasOP}$ and 
$\ref{totalsymbimgcontains10}$ in the Proposition \ref{symboliccharacterizationofOP} shows that whether or not $\phi$ has OP is reflected in 
$\rho_{\phi}(M)$, the $\mathbb{Z}$-total symbolic image of $\phi$.

\subsubsection{Strict order property}
In this part we characterize strict order property in terms of symbolic images of formulas.
Let $\Theta$ be subset of $\mathcal{B}$ consisting of shift orbit of $\bar{1}\bar{0}$. 

\begin{pro} \label{symboliccharacterizationofSOP}
Let $T$ be a theory and $\phi(x,y)$ be a formula. Then the followings are equivalent.

\begin{enumerate}
\item{
The formula $\phi(x,y)$ has strict order property.} \label{phihasSOP}

\item{
There exists some model $M$ of $T$, some instance $U=\phi(M,b)$ for some $b \in M$ and some automorphism $\sigma$ of $M$ such that 
$U \subsetneqq \sigma(U)$.} \label{Ustrictsubsetoftau(U)}

\item{ 
There exists some model $M$ of $T$, some instance $U=\phi(M,b)$ for some $b \in M$ and some automorphism $\sigma$ of $M$ such that $\xi_{\sigma,U}(M)$ consists of $\Theta$ and possibly $\bar{0}$ and $\bar{1}$.} \label{imageconsistsoforbitof10}

\item
{Same as \ref{imageconsistsoforbitof10} with additional property that $\sigma \in W_{b}$.} \label{imageconsistsoforbitof10autinWb}
\end{enumerate}
\end{pro}

\proof
$(\ref{phihasSOP} \Rightarrow \ref{Ustrictsubsetoftau(U)})$ 
By using the Lemma \ref{differentequivalencesofSOP}, there exists some model $M$ of $T$ and some indiscernible sequence $I=(b_i)_{i \in \mathbb{Z}}$ in $M$ such that
$\phi(M,b_i) \subsetneqq \phi(M,b_{i+1})$ for every $i \in \mathbb{Z}$.
We may assume that $M$ is saturated enough. By saturation and indiscerniblity of $I=(b_i)_{i \in \mathbb{Z}}$, there exists some automorphism $\sigma$ of $M$ such that $\sigma^i(b_0)=b_i$ for every $i \in \mathbb{Z}$. Let $U=\phi(M,b_0)$. So we have that $U \subsetneqq \sigma(U)$.

\vspace{1.5mm}

$(\ref{Ustrictsubsetoftau(U)} \Rightarrow \ref{phihasSOP})$
Since $U \subsetneqq \sigma(U)$, we have $\sigma^i(U) \subsetneqq \sigma^{i+1}(U)$ for every $i \in \mathbb{Z}$. So 
$\phi(M,b_i) \subsetneqq \phi(M,b_{i+1})$ for every $i \in \mathbb{Z}$. Hence $\phi$ has SOP.

\vspace{1.5mm}

$(\ref{Ustrictsubsetoftau(U)} \Rightarrow \ref{imageconsistsoforbitof10})$
Since $U \subsetneqq \sigma(U)$ we have $\sigma^{-1}(U) \subsetneqq U$. So there is no appearance of $01$ as a subsequence in $\xi_{\sigma,U}(a)$ for every $a \in M$. So $\xi_{\sigma,U}(M)$ consists of $\Theta$ and possibly $\bar{0}$ and $\bar{1}$.

\vspace{1.5mm}

$(\ref{imageconsistsoforbitof10} \Rightarrow \ref{Ustrictsubsetoftau(U)})$
We show that $U \subsetneqq \sigma(U)$. First we show that $U \subseteq \sigma(U)$. Assume for contradiction that $a \in U \setminus \sigma(U)$. 
Let $I=\xi_{\sigma,U}(a)$. So $I(-1)=0$ and $I(0)=1$. Hence $I$ contains $01$ as a subsequence which is a contradiction. Therefore $U \subseteq \sigma(U)$.
Since $\xi_{\sigma,U}(M)$ contains $\Theta$, there is some $a \in M$ such that $\xi_{\sigma,U}(a)=\bar{1}\bar{0}$. So $a \in U$ but $\sigma(a) \not \in U$ implying $a \not \in \sigma^{-1}(U)$. Hence $\sigma^{-1}(U) \subsetneqq U$. Now we have $U \subsetneqq \sigma(U)$.
$\ \ \square$

\vspace{2mm}

Now we give a characterization for NSOP theories in terms of symbolic representation and SW-closeness.

\begin{pro} \label{nsopssclosed}
A theory $T$ is NSOP if and only if for every (some) $\aleph_1$-saturated model $M$ of $T$, every formula $\phi(x,y)$, every instance $U=\phi(M,b)$ of $\phi$ with parameter $b \in M$ and every $\sigma \in W_{b}$, the symbolic image
$\xi_{\sigma,U}(M)$ is SW-closed in $\mathcal{B}$.
\end{pro}
\proof
$\Rightarrow$) Assume for contradiction that $\xi_{\sigma,U}(M)$ is not SW-closed where $U=\phi(M,b)$ is some instance of some formula in some saturated model $M$ and $\sigma \in W_{b}$.
So, since $\xi_{\sigma,U}(M)$ is shift-invariant,
there are $I,J \in \mathcal{B}$ such that 
$J =w(I)$ and $I \in \xi_{\sigma,U}(M)$ while $J \not \in \xi_{\sigma,U}(M)$.  
Without lose of generality we may assume that $I(0)=1,I(1)=0, J(0)=0$ and $J(1)=1$. 
Let $A:=\{\psi_i, \ i \in \mathbb{Z}\}$ and $B:=\{\chi_i, \ i \in \mathbb{Z}\}$
where $\psi_i=\phi^{I(i)}(\x,\sigma^{-i}(b))$ and $\chi_i=\phi^{J(i)}(\x,\sigma^{-i}(b))$
for every $i \in \mathbb{Z}$.
One sees that $\psi_i=\chi_i$ for every $i \not =0,1$ and so  $A$ and $B$ are different only for $i =0,1$. Also $\psi_0=\neg \chi_1$ and $\psi_1=\neg \chi_0$.
Since $I \in \xi_{\sigma,U}(M)$, $A$ is satisfiable by any $a \in M$ with $I=\xi_{\sigma,U}(a)$. But since  $J \not \in \xi_{\sigma,U}(M)$ and $M$ is saturated, $B$ is not satisfiable and is inconsistent.
So there exists some $n \in \mathbb{N}$ such that for 
$Z:=\bigcap_{i\in [-n,n] \setminus \{0,1\}} \psi^M_i$, 
$X:= \phi(M,b)$ and 
$Y:=\phi(M,\sigma^{-1}(b))$, 
we have that
$X \cap Y^c \cap Z \not= \emptyset$ while $X^c \cap Y \cap Z =\emptyset$.
Therefore $Y \cap Z \subsetneqq X \cap Z$.
Since $\sigma \in W_{b}$, the $\sigma$-orbit of $b$ is an indiscernible sequence. Hence there exists some automorphism $\tau$ such that 
$Z=\tau(Z)$ and $X=\tau(Y)$.
Let $R=Y \cap Z$. So $R$ is a definable set with the property that $R \subsetneqq \tau(R)$. We may assume that $R$ is defined by some instance of some formula $\zeta$. Now by part \ref{Ustrictsubsetoftau(U)} of Proposition \ref{symboliccharacterizationofSOP}, $\zeta$ has SOP which is a contradiction with our assumption that $T$ is NSOP.

\vspace{2mm}

$\Leftarrow$) Assume for contradiction that some formula $\phi(x,y)$ has SOP. 
So by part \ref{imageconsistsoforbitof10autinWb} of Proposition \ref{symboliccharacterizationofSOP},
there exists some model $M$ of $T$, some instance $U=\phi(M,b)$ for some $b \in M$ and some automorphism $\tau \in W_{b}$ such that $\xi_{\tau,U}(M)$ consists of $\Theta$ and possibly $\bar{0}$ and $\bar{1}$. Let 
$I$ be the sequence with $I(i)=1$ for $i < 0$ and $i=1$, and $I(i)=0$ for $i=0$ and $i >1$. Now $I$ is a switching of $\bar{1}\bar{0}$ while $I$ is not in $\xi_{\tau,U}(M)$. This contradicts with SW-closeness of $\xi_{\tau,U}(M)$. 
$\ \ \square$

\subsubsection{Shelah's theorem: OP if and only if IP or SOP}
In this part we conclude Shelah's characterization of OP in term of IP and SOP using the language and materials presented earlier. Note that by comparing part \ref{imageconsistsoforbitof10} of Proposition \ref{symboliccharacterizationofSOP} and
part \ref{orbitclosurecontains10} of Proposition \ref{symboliccharacterizationofOP}, one easily concludes that if $\phi$ has SOP then it has OP.

\begin{thm}
Let $T$ be a theory and $\phi(x,y)$ be a formula with OP. Then either $\phi$ has IP or there exists some formula with SOP.
\end{thm}
\proof
Assume that every formula of $T$ is NSOP.   
Since $\phi(x,y)$ has OP,  
by part part \ref{10imageofsomeawithtauinWb} of Proposition \ref{symboliccharacterizationofOP},
there exists some $b \in M$ and some $\sigma \in W_{b}$ such that $\bar{1}\bar{0} \in \xi_{\sigma,U}(M)$ where $U=\phi(M,b)$. Since $T$ has NSOP, by Proposition \ref{nsopssclosed} $H:=\xi_{\sigma,U}(M)$ is SW-closed.
So $\langle \{\bar{1}\bar{0}\} \rangle \subseteq \langle H \rangle=H$. Hence by Remark \ref{<01>dense}, $H$ is dense in $\mathcal{B}$ with the natural metric. Now by theorem \ref{symboliccharacterizationofNIP} $\phi$ has IP. $\ \ \square$

\bigskip

{\bf Acknowledgement.} 
The author would like to thank L. Newelski and Wroclaw logic group for helpful discussions.
Also the author is indebted to
Institute for Research in Fundamental Sciences, IPM, for support.
This research was in part supported by a grant from IPM (No. 93030035).


\begin{thebibliography}{99}


\bibitem{EntropyinDS} T. Downarowicz, \textit{Entropy in Dynamical Systems}, Cambridge University Press, 2011.


\bibitem{HalmosVonNeuman} P. Halmos, J. Von Neumann, \emph{Operator methods in classical mechanics, II}, Annals of Mathematics. Second Series (Annals of Mathematics) 43 (2): (1942), 332–350.

\bibitem{HPP} E. Hrushovski, Y. Peterzil, and A. Pillay, \emph{Groups, measures and the NIP}, Journal AMS, 21 (2008), 563-596.

\bibitem{HP} E. Hrushovski and A. Pillay, \emph{On NIP and invariant measures}, Journal of European Math. Society, 13 (2011), 1005-1061.


\bibitem{K1} H. J. Keisler, \emph{Measures and forking}, Annals of Pure and Applied Logic 45 (1987), 119-169.


\bibitem{generalizedNormalNumber} F. J. Martinelli, \textit{Construction of generalized normal numbers}, Pacific Journal of Mathematics, 1978.


\bibitem{N} L. Newelski, \emph{Topological dynamic of definable group actions}, Journal of Symbolic Logic, Journal of Symbolic Logic, 74 (2009), pp. 50-72.


\bibitem{Poizatmodeltheory} B. Poizat, \textit{A course in model theory: An introduction to contemporary mathematical logic}, Springer, New York, (2000).


\bibitem{Shelahclassification} S. Shelah, \emph{Classification Theory and the number of non-isomorphic models}, North-Holland, Amsterdam and New York, (1990).



\end{thebibliography}
\end{document}